\documentclass[12pt]{article}

\usepackage{amsbsy}
\usepackage{amssymb}

\begin{document}

\author{Oleg Pikhurko\thanks{Supported by a Research Fellowship, St.\ John's College, Cambridge.}\\
 DPMMS, Centre for Mathematical Sciences\\
 Cambridge University, Cambridge~CB3~0WB, England\\
 E-mail: {\tt O.Pikhurko@dpmms.cam.ac.uk}}

\newcommand{\eqref}[1]{\mbox{\rm(\ref{#1})}}
\newcommand{\pl}{}
\newcommand{\C}[1]{{\protect\cal #1}}
\newcommand{\B}[1]{{\bf #1}}
\newcommand{\I}[1]{{\mathbb #1}}
\renewcommand{\O}[1]{\overline{#1}}
\newcommand{\binom}[2]{{#1\choose #2}}
\newcommand{\e}{\epsilon}
\newcommand{\imply}{\Rightarrow}
\newcommand{\comment}[1]{}
\newcommand{\qed}{\nolinebreak\mbox{\hspace{5 true pt}%
\rule[-0.85 true pt]{3.9 true pt}{8.1 true pt}}}

\newtheorem{theorem}{Theorem}
\newtheorem{lemma}[theorem]{Lemma}
\newtheorem{corollary}[theorem]{Corollary}
\newtheorem{proposition}[theorem]{Proposition}
\newtheorem{conjecture}[theorem]{Conjecture}
\newtheorem{example}[theorem]{Example}
\newtheorem{problem}[theorem]{Problem}

\newcommand{\rcr}[1]{\ref{cr:\pl:#1}}
\newcommand{\claim}[1]{\smallskip\noindent{\bf Claim #1} }

\newcommand{\V}[1]{{\boldsymbol #1}} 
\newcommand{\VV}[1]{{#1}}
\newcommand{\diam}{\mbox{\rm diam}}
\newcommand{\Times}[1]{{}^{{(}{\times}#1{)}}}

\title{Borsuk's Conjecture Fails in\\ Dimensions $321$ and $322$}
\maketitle

Let $f(n)$ be the smallest $f$ such that any bounded set in $\I R^n$
can be partitioned into at most $f$ sets of smaller diameter.  The
famous Borsuk's conjecture~\cite{borsuk:33} that $f(n)=n+1$ for any
$n\ge 1$ has been spectacularly disproved by Kahn and
Kalai~\cite{kahn+kalai:93}. However, the counterexamples
in~\cite{kahn+kalai:93} all have very large dimension. If we define
 $$
 n_0=\min\{n\in\I N\mid f(n)>n+1\},
 $$
 then the proof of Kahn and Kalai gives $n_0\le 1325$.

On the other hand we know only that $n_0\ge 4$
(Perkal~\cite{perkal:47}; Eggleston~\cite{eggleston:55}). It is of
interest where $n_0$ lies. The upper bound on $n_0$ was improved to
$n_0\le946$ (Nilli~\cite{nilli:94}), $n_0\le561$
(Raigorodski~\cite{raigorodski:97}), $n_0\le 560$ (Wei\ss
bach~\cite{weissbach:00}), and $n_0\le 323$
(Hinrichs~\cite{hinrichs:02}). In fact, we know that $f(n)>n+1$ for
all $n\ge 323$.

Here we show that $n_0\le 321$. 

\begin{theorem}\label{th:\pl:1} $f(321)\ge 333$. Thus, Borsuk's conjecture fails in
all dimensions $n\ge 321$.\end{theorem}

Let us first recall Hinrichs' construction which utilises
$M=\Lambda_{24}\cap \Omega^{24}$, the set of unit-length elements
in the Leech lattice $\Lambda_{24}$. Namely, choose an orthonormal
basis
 $$\left((\V e_i)_{i=1}^{24}, (\V f_i)_{i=1}^{24}, (\V g_{i,j})_{1\le
i<j\le 24}\right)$$
 in $R^{324}$, define $\Phi:\I R^{24}\to \I R^{324}$ by
 $$
 \Phi(\VV x_1,\dots,\VV x_{24})=\frac2{\sqrt5}\, \sum_{i=1}^{24} \VV
x_i^2\V e_i + \frac1{\sqrt5}\, \sum_{i=1}^{24} \VV x_i\V f_i +
\frac{2\sqrt2}{\sqrt5}\sum_{1\le i<j\le 24} \VV x_i\VV x_j\V g_{i,j},
 $$
 and consider $\Phi(M)$.
Hinrichs~\cite[Proposition~3(iii)]{hinrichs:02} proved that any subset
of $\Phi(M)$ of smaller diameter has at most $350$ elements.

To establish Theorem~\ref{th:\pl:1} we show that there is $N\subset M$ such
that $|N|\ge 116424$ and $\Phi(N)$ lies within a $321$-dimensional
affine subspace of $\I R^{324}$. Then
by~\cite[Proposition~3(iii)]{hinrichs:02} we have
 $$\diam(\Phi(N))=\diam(\Phi(M)).$$
 Applying~\cite[Proposition~3(iii)]{hinrichs:02} again, we conclude
that we need at least
 \begin{equation}\label{eq:\pl:321}
 |\Phi(N)|/350=|N|/350\ge 116424/350> 332,
 \end{equation}
 parts of smaller diameter to partition $\Phi(N)$, which implies
the desired inequality $f(321)\ge 333$.

In order to prove the existence of $N$ we need the following explicit
description of $M$, taken from Conway and Sloane~\cite[\S 11~of
Chapter~4]{conway+sloane:splg}. Namely, $M$ contains
 \begin{itemize}
 \item $97152$ points of the form
$\frac1{4\sqrt2}(\pm2\Times8,0\Times{16})$,
 \item $98304$ points of the form $\frac1{4\sqrt2}(\pm
3,\pm1\Times{23})$,
 \item $1104$ points of the form
$\frac1{4\sqrt2}(\pm4\Times2,0\Times{22})$,
 \end{itemize}
 where the actual signs and positions of coordinates are not relevant
for our purposes. ($a\Times{k}$ denotes $k$ copies of~$a$.)

Define the bipartite graph $G$ with parts $M$ and
$\binom{\{1,\dots,24\}}3$ so that $(\VV x_1,\dots,\VV x_{24})$ and
$\{k,l,m\}$ are connected if $|\VV x_{k}|=|\VV x_{l}|=|\VV
x_{m}|$. Clearly,
 $$
 \textstyle  e(G)=97152\times\left(\binom83+\binom{16}3\right) +
98304\times \binom{23}3 + 1104\times \binom{22}3=235642176.
 $$
 Hence, some set $\{k,l,m\}$ receives at least
$e(G)/\binom{24}3=116424$ edges. Let $N\subset M$ consist of its
neighbours. We have
 $$
 \Phi(N)\subset\left\{\V y\in\I R^{324}\mid \V e_{k}\cdot \V y=\V
e_{l}\cdot \V y=\V e_{m}\cdot \V y\mbox{ and } {\textstyle
\sum_{i=1}^{24} \V e_i\cdot \V y=1}\right\},
 $$
 the latter set being a $321$-dimensional affine subspace of $\I
R^{324}$. Thus $N$ has all the required properties.

\medskip\noindent{\bf Remark}  Similarly, one can find a set $K\subset M$
of size $143136$ such that for some $1\le k<l\le 24$ and for any $(\VV
x_1,\dots,\VV x_{24})\in K$ we have $|\VV x_{k}|=|\VV x_{l}|$. Then
$\Phi(K)$ is `$322$-dimensional' and we have $f(322)\ge \lceil\,
|K|/350\,\rceil \ge 409$.\medskip

\noindent{\bf Remark} Theorem~\ref{th:\pl:1} has been independently
discovered by Hinrichs and Richter. Moreover, they report to have
proved $n_0\le 298$ by showing that the set
 $$
 L=\{(x_1,\dots,x_{24})\in M\mid x_1=x_2\}.
 $$
 which lies in a $298$-dimensional affine subspace cannot be
partitioned into $299$ parts of smaller diameter. However, the proof
of the latter claim (being currently writen) seems to be long and
complicated.

\small
\enlargethispage*{20pt}

\providecommand{\bysame}{\leavevmode\hbox to3em{\hrulefill}\thinspace}
\providecommand{\MR}{\relax\ifhmode\unskip\space\fi MR }
\providecommand{\MRhref}[2]{%
  \href{http://www.ams.org/mathscinet-getitem?mr=#1}{#2}
}
\providecommand{\href}[2]{#2}

\end{document}